\documentclass[reqno,11pt]{amsart}
\usepackage{graphicx}
\usepackage{amsmath,amssymb}
\newtheorem{theorem}{Theorem}[section]

\newtheorem{lemma}[theorem]{Lemma}
\newtheorem{proposition}[theorem]{Proposition}
\theoremstyle{definition}
\newtheorem{definition}[theorem]{Definition}
\theoremstyle{remark}

\numberwithin{equation}{section}

\begin{document}

\title{Free Jump Dynamics in Continuum}%
\author{Joanna Bara\'nska}%
\address{Institute of Mathematics, Maria Curie-Sk{\l}odowska University, 20-031 Lublin, Poland}%
\email{asia13p@wp.pl}%

\author{ Yuri Kozitsky}
\address{Institute of Mathematics, Maria Curie-Sk{\l}odowska University, 20-031 Lublin, Poland}%
\email{jkozi@hektor.umcs.lublin.pl}%

\thanks{This work was supported in part by the DFG through SFB 701: ``Spektrale
Strukturen und Topologische Methoden in der Mathematik" and by the European Commission under the project
STREVCOMS PIRSES-2013-612669.}%

\begin{abstract}
The evolution is described of an infinite system of hopping point  particles in $\mathbb{R}^d$. The states of the system are
probability measures on the space of configurations of particles. Under the condition that the initial state $\mu_0$ has correlation functions
of all orders which are: (a) $k_{\mu_0}^{(n)} \in L^\infty ((\mathbb{R}^d)^n)$ (essentially bounded); (b) $\|k_{\mu_0}^{(n)}\|_{ L^\infty ((\mathbb{R}^d)^n)} \leq C^n$, $n\in \mathbb{N}$ (sub-Poissonian), the evolution $\mu_0 \mapsto \mu_t$, $t>0$, is obtained as a continuously differentiable map $k_{\mu_0} \mapsto k_t$, $k_t =(k_t^{(n)})_{n\in \mathbb{N}}$, in the space of essentially bounded sub-Poissonian functions.  In particular, it is proved that $k_t$ solves the corresponding evolution equation, and that for each $t>0$ it is the correlation function of a unique state $\mu_t$.

\end{abstract}
\maketitle
\section{Introduction}
\label{S1}

In this paper, we study the dynamics of an infinite system of point
particles $x\in \mathbb{R}^d$, $d\geq 1$. States of the
system are discrete subsets of $\mathbb{R}^d$ -- {\em
configurations}, which constitute the set
\begin{equation}
  \label{1}
\Gamma = \{ \gamma \subset \mathbb{R}^d: |\gamma \cap \Lambda| <
\infty \  {\rm for}  \  {\rm any}  \ {\rm compact}  \ \Lambda
\subset \mathbb{R}^d \},
\end{equation}
where $|\cdot|$ stands for cardinality. Note that  $\Gamma$ contains
also finite configurations, including the empty one. The set
$\Gamma$ can be completely and separably metrized, see
\cite{Albev,Tobi}, and thus equipped with the corresponding Borel
$\sigma$-field $\mathcal{B}(\Gamma)$. The elements of $\Gamma$ are
considered as point states of the system in the sense that, for  a suitable
function $F:\Gamma\to \mathbb{R}$, the number $F(\gamma)$ is treated
as the value of {\em observable} $F$ in state $\gamma$. Along with
point states $\gamma\in \Gamma$ one  employs states determined
by probability measures on $\mathcal{B}(\Gamma)$. In this case, the
corresponding value is the integral
\begin{equation*}
 \langle \! \langle F, \mu \rangle \! \rangle := \int_{\Gamma} F d \mu,
\end{equation*}
and the system's dynamics are described as  maps $\mu_0 \mapsto
\mu_t$, $t>0$.

In the Markov approach, the map
$\mu_0 \mapsto \mu_t$ is obtained from the Fokker-Planck equation
\begin{equation}
  \label{2}
  \frac{d}{dt} \mu_t = L^\mu \mu_t, \qquad \mu_t|_{t=0} = \mu_0, \quad t>0,
\end{equation}
in which `generator' $L^\mu$ specifies the model. In order to solve
(\ref{2}) in the set of probability measures one has to introduce an
appropriate mathematical setting, e.g., a Banach space of signed
measures, and then to define $L^\mu$ as a linear operator in this
space. However, for infinite systems such a direct way is rather
impossible. By the duality
\begin{equation}
 \label{R}
\langle \! \langle F_0, \mu_t \rangle \! \rangle = \langle \!
\langle F_t, \mu_0 \rangle \! \rangle, \qquad t>0,
\end{equation}
the observed evolution $\langle \! \langle F, \mu_0
\rangle \!  \rangle \mapsto \langle \! \langle F, \mu_t \rangle \!
\rangle$
can also be considered as the evolution $\langle \! \langle F_0, \mu \rangle \!  \rangle
\mapsto \langle \! \langle F_t, \mu \rangle \!  \rangle$ obtained from the Kolmogorov equation
\begin{equation}
 \label{R2}
\frac{d}{d t} F_t = L F_t , \qquad F_t|_{t=0} = F_0, \qquad t>0,
\end{equation}
where $L$ and $L^\mu$ are dual in the sense of (\ref{R}). Thus, also
`generator' $L$ specifies the model. Various types of such
generators are discussed in \cite{Dima2}. In this paper, we consider the model specified by
\begin{equation}
  \label{R20}
  (LF)(\gamma) = \sum_{x\in \gamma} \int_{\mathbb{R}^d} a(x-y) \left[F(\gamma \setminus x \cup y) - F(\gamma) \right]d y.
\end{equation}
Here and in the sequel in the corresponding context,  $x\in
\mathbb{R}^d$ is also treated as a single-point configuration
$\{x\}$. The {\em jump} kernel $a(x) = a(-x)\geq 0$ is
supposed to satisfy the condition
\begin{equation}
  \label{a}
   \int_{\mathbb{R}^d} a (x) d x =: \alpha < +\infty.
\end{equation}
The `generator' in (\ref{R20}) describes free jumps of the elements of
configurations. In models where jumps are not free, the kernel $a$
depends also on $\gamma$, see \cite{Berns,Dima2}.

Similarly as above, to solve (\ref{R2}) one should define $L$ as a
linear operator in an appropriate Banach space of functions, which
can also be problematic as the sum in (\ref{R20}) typically runs
over an infinite set. One of the possibilities here is to construct
a Markov process $\gamma_t$ with state space $\Gamma$, which starts
from a fixed configuration $\gamma_0\in \Gamma$. Then
$\mu^{\gamma_0}_t$ -- the law of $\gamma_t$, solves (\ref{2}) with
$\mu_0=\delta_{\gamma_0}$ (the Dirac measure). However, since the
evolution of the model (\ref{R20}) includes simultaneous jumps of an infinite
number of points, there can exist $\gamma_0$ such that with
probability one at some $t>0$ infinitely many points appear in
a bounded $\Lambda$, see the corresponding discussion in
\cite{icmp}. The reason for this is that the
configuration space (\ref{1}) appears to be too big and cannot serve
as a state-space for the corresponding process. In \cite{icmp},
under a more restrictive condition than that in (\ref{a}) (see eq.
(39) in that paper), the Markov process corresponding to (\ref{R20})
was constructed for $\gamma_0$, and hence all $\gamma_t$, $t>0$, lying
in a certain proper measurable subset $\Theta \subset\Gamma$. By
this result, the evolution $\mu_0 \mapsto \mu_t$ corresponding to
(\ref{2}) with $\mu_0 (\Theta)=1$ can be obtained by the formula
\[
\mu_t (\cdot) = \int_{\Gamma} \mu^{\gamma}_t (\cdot) \mu_0(
d\gamma),
\]
which guarantees also that $\mu_t (\Theta) =1$ for all $t>0$.

There exists another approach to solving (\ref{2}) in which instead
of restricting the set of configurations where the process takes
its values one restricts the set of initial measures $\mu_0$. This
restriction amounts to imposing a condition, formulated in terms of
the so called {\em correlation measures}. For $n\in \mathbb{N}$ and
a probability measure $\mu$ on $\mathcal{B}(\Gamma)$, the $n$-th order correlation measure
$\chi^{(n)}_\mu$ is related to $\mu$ by the following formula
\begin{eqnarray}
  \label{3}
& & \int_{\Gamma} \bigg{(} \sum_{\{x_1, \dots , x_n\}\subset \gamma} G^{(n)}(x_1 , \dots , x_n) \bigg{)} \mu( d \gamma) \\[.2cm]
& & \qquad = \frac{1}{n!} \int_{(\mathbb{R}^d)^n} G^{(n)}(x_1 , \dots , x_n) \chi^{(n)}_\mu(dx_1 , \dots , d x_n), \nonumber
\end{eqnarray}
which ought to hold for all bounded, compactly supported  measurable functions
$G^{(n)}:(\mathbb{R}^d)^n \to \mathbb{R}$, see, e.g.,
\cite{Dima2,Tobi}. Now the mentioned condition is that, for each
$n\in \mathbb{N}$,  $\chi^{(n)}_\mu$ is absolutely continuous with
respect to the corresponding Lebesgue measure with Radon-Nikodym
derivative $k^{(n)}_\mu\in L^{\infty}((\mathbb{R}^d)^n)$. Clearly,
this condition excludes Dirac measures, and hence the possibility to
solve (\ref{2}) by means of stochastic processes. Instead one can
drop the mentioned above restrictions on the jump kernel $a$ and on
the support of $\mu_t$. In the present paper, we follow this way.

Until this time, there have been published only two papers
\cite{Berns,icmp} dealing with jump models on $\Gamma$. In \cite{Berns},
the approach based on essentially bounded correlation functions is
applied to the model in which the jump kernel depend also on the
configuration $\gamma \setminus x$ in a specific way. In that paper,
by means of the so called {\em Ovcyannikov} method \cite{Trev}, the
evolution $\mu_0 \mapsto \mu_t$ is constructed for $t\in [0,T)$,
with $T< \infty$ dependent on the kernel $a$. In the present paper, we
consider the free case (no dependence of $a$ on $\gamma$),
which allows us to employ semigroup methods and obtain the
evolution $\mu_0 \mapsto \mu_t$ for all $t\geq 0$.

\section{Basic Notions and the Result}
\label{S2}

The main idea of the approach used in this paper is to obtain the
evolution of states from the evolution of their correlation
functions. This includes the following steps: (a) passing from
problem (\ref{2}) to the corresponding problem for correlation
functions; (b) obtaining $k^{(n)}_{\mu_0}=: k_0^{(n)} \mapsto
k_t^{(n)}$; (c) proving that, for each $t>0$, there exists a unique
$\mu_t$ such that $k_t^{(n)}$ is its correlation function for all
$n\in \mathbb{N}$. We perform this
steps in subsection \ref{SS22}. In subsection \ref{SS21}, we present some details
of the method. Further information on the methods used in this work can be found
in \cite{Berns,FKK,Dima,Dima2,Tobi}.

\subsection{Configuration spaces and correlation functions}
\label{SS21}

By $\mathcal{B}(\mathbb{R}^d)$ and $\mathcal{B}_{\rm
b}(\mathbb{R}^d)$ we denote the set of all Borel and all bounded
Borel subsets of $\mathbb{R}^d$, respectively. The set of
configurations $\Gamma$ defined in (\ref{1}) is equipped with the
vague topology -- the weakest topology which makes the maps
\[
\Gamma \ni \gamma \mapsto \sum_{x \in \gamma} f(x) \in \mathbb{R}
\]
continuous for all compactly supported continuous functions
$f:\mathbb{R}\to \mathbb{R}$. This topology can be completely and
separably metrized, that turns $\Gamma$ into a Polish spaces, see
\cite{Albev,Tobi}. By $\mathcal{B}(\Gamma)$ and
$\mathcal{P}(\Gamma)$ we denote the Borel $\sigma$-field of subsets
of $\Gamma$ and the set of all probability measures on
$\mathcal{B}(\Gamma)$, respectively.

The set of finite configurations
\begin{equation}
  \label{6}
  \Gamma_{0} = \bigsqcup_{n\in \mathbb{N}_0} \Gamma^{(n)} ,
\end{equation}
is the disjoint union of the sets of $n$-particle configurations:
\begin{equation*}
\Gamma^{(0)} = \{ \emptyset\}, \qquad \Gamma^{(n)} = \{\gamma \in \Gamma: |\gamma| = n \}, \ \ n\in \mathbb{N}.
\end{equation*}
For $n\geq 2$, $\Gamma^{(n)}$ can be identified with the symmetrization of the set
\begin{equation}
  \label{7a}
\left\{(x_1, \dots , x_n)\in \bigl(\mathbb{R}^{d}\bigr)^n: x_i \neq x_j, \ {\rm for} \ i\neq j\right\},
\end{equation}
which allows one to introduce the corresponding (Euclidean) topology
on $\Gamma^{(n)}$. Then by (\ref{6}) one defines also the topology
on the whole $\Gamma_0$: $A\subset \Gamma_0$ is said to be open if its
intersection with each $\Gamma^{(n)}$ is open. This topology differs
from that induced on $\Gamma_0$ by the vague topology of $\Gamma$.
At the same time, as a  set $\Gamma_0$  is in
$\mathcal{B}(\Gamma)$. Thus, a function $G :\Gamma_0 \to \mathbb{R}$
is measurable as a function on $\Gamma$ if and only if its
restrictions to each $\Gamma^{(n)}$ are Borel functions. Clearly,
these restrictions fully determine $G$. In view of (\ref{7a}), the
restriction of $G$ to $\Gamma^{(n)}$ can be extended to  a symmetric function
$G^{(n)}: (\mathbb{R}^d)^n \to \mathbb{R}$, $n\in \mathbb{N}$, such
that
\begin{equation}
  \label{8}
 G(\gamma) = G^{(n)} (x_1, \dots , x_n), \qquad {\rm for} \ \  \gamma = \{x_1 , \dots , x_n\}.
\end{equation}
It is convenient to complement (\ref{8})  by putting
$G(\emptyset)=G^{(0)} \in \mathbb{R}$.
\begin{definition}
  \label{bbsdf}
A measurable function $G:\Gamma_0 \to \mathbb{R}$ is said to have {\em bounded support}
if the following holds: (a)  there exists $N\in \mathbb{N}_0$ such that $G^{(n)} \equiv 0$ for all $n> N$; (b) there exists
$\Lambda \in \mathcal{B}_{\rm b}(\mathbb{R}^d)$ such that, for all $n\in \mathbb{N}$, $G^{(n)} (x_1 , \dots , x_n) = 0$  whenever $x_j \in  \mathbb{R}^d \setminus \Lambda$ for some $j=1 , \dots , n$. By $B_{\rm bs} (\Gamma_0)$ we denote the set of all such functions.
\end{definition}
Let all $G^{(n)}$, $n\in \mathbb{N}$, be
bounded Borel functions and $G$ be related to $G^{(n)}$ by (\ref{8}).
For such $G$, we then write
\begin{eqnarray}
  \label{9}
\int_{\Gamma_0}  G(\gamma) \lambda (d \gamma) = G^{(0)} + \sum_{n=1}^\infty \frac{1}{n!} \int_{(\mathbb{R}^d)^n} G^{(n)} (x_1 , \dots x_n) d x_1 \cdots d x_n .
\end{eqnarray}
This expression determines a
$\sigma$-finite measure $\lambda$ on $\Gamma_0$, called the {\em
Lebesgue-Poisson} measure. Then the formula in (\ref{3}) can be
written in the following way
\begin{equation}
  \label{10}
\int_{\Gamma} \left( \sum_{\eta \Subset \gamma} G(\eta)\right) \mu(d \gamma) = \int_{\Gamma_0} G(\eta) k_\mu (\eta) \lambda (d \eta),
\end{equation}
where the sum on the left-hand side runs over all finite
sub-configurations of $\gamma$. Like in (\ref{8}), $k_\mu : \Gamma_0
\to \mathbb{R}$ is determined by its restrictions $k^{(n)}_\mu$.
Note that $k^{(0)}_\mu \equiv 1$ for all $\mu\in
\mathcal{P}(\Gamma)$.

Prototype examples of
measures with property $k^{(n)}_\mu \in L^\infty ((\mathbb{R}^d)^n)$  are the Poisson measures
$\pi_\varrho$ for which
\begin{equation}
  \label{5}
k^{(n)}_{\pi_\varrho} (x_1, \dots , x_n) = \prod_{j=1}^n
\varrho(x_j), \qquad n\in \mathbb{N}.
\end{equation}
Here  $\varrho\in L^\infty (\mathbb{R}^d)$, and the case of constant $\varrho
\equiv \varkappa >0$ corresponds to the homogeneous Poisson measure.
Along with the spatial properties of correlation functions, it
is important to know how do they depend on $n$. Having in mind
(\ref{5}) we say that a given $\mu\in \mathcal{P}(\Gamma)$ is {\em
sub-Poissonian} if its correlation functions are such that
\begin{equation}
  \label{11}
 k_\mu^{(n)} (x_ 1 , \dots , x_n) \leq C^n ,
\end{equation}
holding for some $C>0$, all $n\in
\mathbb{N}$, and Lebesgue-almost all $(x_1, \dots, x_n)\in(\mathbb{R}^d)^n$.
Then a  state with property (\ref{11}) is similar to the Poisson states $\pi_\varrho$ in which
the particles are independently scattered over $\mathbb{R}^d$. At the same time, the increase of
$k^{(n)}_\mu$ as $n!$ corresponds to the appearance of {\it clusters} in state $\mu$.
For the so called continuum contact model, it is
known \cite{Dima} that, for any $t>0$,
\begin{equation*}
{\rm const}\cdot n! c^n_t \leq
 k^{(n)}_t (x_1, \dots, x_n) \leq {\rm const}\cdot n! C^n_t,
\end{equation*}
where the left-hand inequality holds if all $x_i$ belong to a ball
of small enough radius.

Recall that by $B_{\rm bs} (\Gamma_0)$ we denote the set of all
$G:\Gamma_0 \to \mathbb{R}$ which have bounded support, see Definition \ref{bbsdf}. For each
such $G$ and $\gamma \in \Gamma$, the expression
\begin{equation}
  \label{bsa}
(K G )(\gamma) := \sum_{\eta \Subset \gamma} G(\eta)
\end{equation}
is well-defined as the sum has finitely many terms only. Note that $KG$ is $\mathcal{B}(\Gamma_0)$-measurable for
each $G\in  B_{\rm bs} (\Gamma_0)$. Indeed, given $G\in B_{\rm bs}(\Gamma_0)$, let $N$ and $\Lambda$ be as in Definition \ref{bbsdf}.
Then
\[
(K G )(\gamma) = G^{(0)} + \sum_{x\in \gamma\cap \Lambda} G^{(1)} (x) + \cdots  + \sum_{\{x_1, \dots x_N\}\subset \gamma\cap \Lambda} G^{(N)} (x_1, \dots, x_N).
\]
The inverse of (\ref{bsa}) has the form
\begin{equation}
  \label{bsb}
(K^{-1} F )(\eta) := \sum_{\xi \subset \eta} (-1)^{|\eta\setminus \xi|}F(\xi).
\end{equation}
As was shown in \cite{Tobi}, $K$ and $K^{-1}$ are linear isomorphisms between $B_{\rm bs}(\Gamma_0)$ and the set of {\em cylinder} functions
$F:\gamma \to \mathbb{R}$.

Let us now turn to the following analog of the classical moment
problem: given a function $k:\Gamma_0 \to \mathbb{R}$, which
properties of $k$ could guarantee that there exists $\mu \in
\mathcal{P}(\Gamma)$ such that $k = k_\mu$? The answer to this
question is given by the following statement, see Theorems 6.1, 6.2
and Remark 6.3 in \cite{Tobi}, in which
\begin{equation}
  \label{bs}
B^{+}_{\rm bs} (\Gamma_0) :=\{ G \in B_{\rm bs} (\Gamma_0): (KG)(\gamma) \geq 0 , \quad \gamma \in \Gamma\} .
\end{equation}
\begin{proposition}
  \label{1pn}
 Let  $k:\Gamma_0 \to \mathbb{R}$ be such that: (a) $k^{(0)} \equiv 1$ and for each $G\in B^{+}_{\rm bs} (\Gamma_0)$ the following holds
 \begin{equation}
   \label{cf}
 \langle \! \langle G, k \rangle \! \rangle := \int_{\Gamma_0} G(\eta) k(\eta) \lambda (d \eta) \geq 0;
 \end{equation}
(b) there exists $C>0$ such that each $k^{(n)}$, $n\in \mathbb{N}$,
satisfies (\ref{11}). Then there exists a unique $\mu\in
\mathcal{P}(\Gamma)$ such that $k$ is its correlation function,
i.e., it is the Radon-Nikodym derivative of the corresponding
correlation measure $\chi_\mu$.
\end{proposition}
Note that $B^+_{\rm bs}(\Gamma_0)$ contains not only positive functions, cf. (\ref{bsb}). That is, the positivity of $k$ as in (\ref{cf}), which readily follows from
(\ref{10}), (\ref{bsa}), and (\ref{bs}), is a stronger property than the usual positivity.

\subsection{The evolution equation}
\label{SS22}
If we rewrite (\ref{10}) in the form, cf. (\ref{bsa}),
\[
\int_{\Gamma} (K G) (\gamma) \mu(d\gamma) = \int_{\Gamma_0} G(\eta) k_\mu (\eta) \lambda (d \eta),
\]
then the action of $L$ on $F$ in (\ref{R20}) can be transferred to $G$, and then to $k_\mu$, as follows
\begin{eqnarray*}
  \int_{\Gamma} [L(K G)] (\gamma) \mu(d\gamma)& = &
  \int_{\Gamma} [K (\widehat{L} G)] (\gamma) \mu(d\gamma)\\[.2cm] & = &
  \int_{\Gamma_0}( \widehat{L} G)(\eta) k_\mu (\eta) \lambda (d \eta) \\[.2cm]  & = &
\int_{\Gamma_0}G(\eta)(L^\Delta) k_\mu (\eta) \lambda (d \eta).
\end{eqnarray*}
Thus, by (\ref{bsa}) and (\ref{bsb}), we see that
\[
\widehat{L} = K^{-1} L K,
\]
and that $L^\Delta$ is the adjoint of $\widehat{L}$ with respect to the pairing in (\ref{cf}).
Then the problem in (\ref{R2}) is being transformed into
the following one
\begin{equation}
  \label{12}
\frac{d}{d t} k_t = L^\Delta k_t, \qquad k_t|_{t=0} = k_0 = k_{\mu_0}.
\end{equation}
For a more general model of jumps in $\mathbb{R}^d$, the calculations of $\widehat{L}$ and $L^\Delta$ were performed in \cite[Section 4]{Dima2}.
The peculiarity of our simple case is that the action of both these `operators' is the same, and, in fact, the same as that in (\ref{R20}). That is,
\begin{gather}
  \label{l}
L^\Delta = A^\Delta + B^\Delta, \\[.2cm]
(A^\Delta k) (\eta) = - \alpha |\eta| k(\eta), \nonumber \\[.2cm]
 (B^\Delta k) (\eta) = \sum_{x\in \eta} \int_{\mathbb{R}^d} a (x-y) k(\eta \setminus x \cup y) d y, \qquad \eta \in \Gamma_0. \nonumber
\end{gather}
So far, these are only informal expressions, like the one in
(\ref{R20}), and our aim now is to define the corresponding linear
operator in a Banach space, where we then solve (\ref{12}).
 Note that in contrast to (\ref{R20}) the sum
in the last line of (\ref{l}) is finite, which shows the
advantage of the approach we follow in this work.
Note also that the action of both $A^\Delta$ and $B^\Delta$ on the
elements of $B_{\rm bs}(\Gamma_{0})$ is well-defined.

In view of our basic assumption $k^{(n)}_\mu\in
L^{\infty}((\mathbb{R}^d)^n)$, the space in question is defined as
follows: for a $\vartheta \in \mathbb{R}$ and a function $u:
\Gamma_0 \to \mathbb{R}$ such that $u^{(n)} \in
L^{\infty}((\mathbb{R}^d)^n)$ for all $n\in \mathbb{N}$, we set
\[
\|u\|_{\vartheta, \infty} = \sup_{n \in \mathbb{N}_0} e^{\vartheta
n} \|u^{(n)}\|_{L^{\infty}((\mathbb{R}^d)^n)},
\]
that can also be written in the form
\begin{equation}
  \label{q1}
\|u\|_{\vartheta, \infty} = \mathrm{ess\,sup}_{\eta \in \Gamma_0} |u(\eta)| \exp( \vartheta |\eta|).
\end{equation}
By means of (\ref{q1}) we then define
\begin{equation*}
  \mathcal{K}_{\vartheta} = \{ u: \Gamma_0 \to \mathbb{R}: \|u\|_{\vartheta, \infty} < \infty\},
\end{equation*}
which is a real Banach space with the standard point-wise linear
operations. Along with (\ref{q1}) we introduce the following norm of
$v: \Gamma_0 \to \mathbb{R}$, cf. (\ref{9}),
\begin{equation}
  \label{q2a}
  \|v \|_{\vartheta, 1} = \int_{\Gamma_0} |v(\eta)| \exp(- \vartheta |\eta|) \lambda (d\eta),
\end{equation}
and then
\begin{equation}
  \label{q2b}
  \mathcal{G}_{\vartheta} = \{ v: \Gamma_0 \to \mathbb{R}: \|v\|_{\vartheta, 1} < \infty\},
\end{equation}
which is also a real Banach space. Note that $\mathcal{K}_{\vartheta}$ is the topological dual to (\ref{q2b}) with the pairing
(\ref{cf}).
To proceed further we need the
following formula, see, e.g., \cite[Lemma 2.1]{Dima},
\begin{eqnarray}
  \label{21}
& & \int_{\Gamma_0} \left( \int_{\mathbb{R}^d} f(y,\eta) g ( \eta) d y \right)\lambda
(d \gamma)\\[.2cm] & & \qquad = \int_{\Gamma_0} \left(\sum_{x \in \eta} f(x,\eta\setminus x)
g ( \eta\setminus x) \right)\lambda( d\eta), \nonumber
\end{eqnarray}
which holds for all appropriate functions $f: \mathbb{R}^d \times
\Gamma_0\to \mathbb{R}$ and $g: \Gamma_0\to \mathbb{R}$.
A more general form of this identity is
\begin{equation}
  \label{21a}
 \int_{\Gamma_0} \left(\sum_{\xi \subset \eta} f(\xi) \right) g(\eta) \lambda ( d\eta) = \int_{\Gamma_0} \int_{\Gamma_0} f(\xi) g (\eta \cup \xi) \lambda (\ d \xi) \lambda ( d \eta).
\end{equation}
Consider
\begin{equation}
  \label{q3}
 \widehat{ \mathcal{G}}_{\vartheta} = \{ v \in \mathcal{G}_\vartheta: A^\Delta v \in  \mathcal{G}_\vartheta\},
\end{equation}
which is a dense linear subset of $\mathcal{G}_{\vartheta}$. Now let
us define in $\mathcal{G}_\vartheta$ linear operators
$\widehat{A}_\vartheta$ and $\widehat{B}_\vartheta$ by the
expressions for $A^\Delta$ and $B^\Delta$, respectively, given in
(\ref{l}). Namely, the action of $\widehat{A}_\vartheta$ and
$\widehat{B}_\vartheta$ on the elements of
$\widehat{\mathcal{G}}_\vartheta$ is given by the right-hand sides
of the second and third expressions in (\ref{l}), respectively.
Clearly $\widehat{A}_\vartheta$  with domain $ \widehat{
\mathcal{G}}_{\vartheta}$ is a closed linear operator in
$\mathcal{G}_{\vartheta}$. Moreover, by (\ref{21}) we get
\begin{eqnarray}
  \label{q4}
 \|{B}^\Delta v\|_{\vartheta,1} &\leq & \int_{\Gamma_0} \int_{\mathbb{R}^d} \sum_{x\in \eta} a(x-y) |v(\eta \setminus x \cup y) | \exp(-\vartheta |\eta|) \lambda ( d\eta)dy \nonumber \\[.2cm]
 & = &  \int_{\Gamma_0} \int_{\mathbb{R}^d} \sum_{x\in \eta} a(x-y) |v(\eta ) | \exp(-\vartheta |\eta|) \lambda ( d\eta)dy  \\[.2cm]
 & = & \alpha \int_{\Gamma_0} |\eta| |v(\eta ) | \exp(-\vartheta |\eta|) \lambda ( d\eta) = \|{A}^\Delta v\|_{\vartheta,1}, \nonumber
\end{eqnarray}
 where we have taken into account (\ref{a}) and  twice used
 (\ref{21}). Thus, $(\widehat{B}_\vartheta, \widehat{\mathcal{G}}_\vartheta)$
 is also well-defined, which allows us to define the sum
$\widehat{L}_\vartheta = \widehat{A}_\vartheta +
\widehat{B}_\vartheta$ with domain $ \widehat{
\mathcal{G}}_{\vartheta}$.  It can readily be shown that
$(\widehat{L}_\vartheta, \widehat{ \mathcal{G}}_{\vartheta})$ is
closed, and its action on the the elements of $\widehat{
\mathcal{G}}_{\vartheta}$ is
\begin{equation}
 \label{q5}
 (\widehat{L}_\vartheta v)(\eta) = - \alpha |\eta| v(\eta) + \int_{\mathbb{R}^d} \sum_{x\in \eta} a(x-y) v(\eta \setminus x \cup y) d
 y.
 \end{equation}
By (\ref{21}), for $v\in \widehat{ \mathcal{G}}_{\vartheta}$ and
$u\in B_{\rm bs}(\Gamma_{0})$, we get
\begin{equation}
  \label{q3a}
  \langle \! \langle \widehat{A}_\vartheta v ,u  \rangle \! \rangle =  \langle \! \langle v, A^\Delta u\rangle \! \rangle,
  \qquad \langle \! \langle \widehat{B}_\vartheta v ,u  \rangle \! \rangle =  \langle \! \langle v, B^\Delta u\rangle \!
  \rangle,
\end{equation}
that is, the action of the dual operators is again given by the same
expressions (\ref{l}). We use this fact to define the adjoint
operator $\widehat{L}^*_\vartheta$. Its domain is
\begin{equation}
  \label{q6}
  \mathcal{D}_\vartheta = \{ u \in \mathcal{K}_\vartheta: \forall v \in  \widehat{ \mathcal{G}}_\vartheta \  \exists w \in \mathcal{K}_\vartheta \
  \langle \! \langle \widehat{L}_\vartheta v ,u  \rangle \! \rangle = \langle \! \langle v ,w  \rangle \! \rangle \}.
\end{equation}
Then the action of $\widehat{L}^*_\vartheta$ on the elements of
$\mathcal{D}_\vartheta \subset \mathcal{K}_\vartheta$ is again given
by (\ref{l}) or by the right-hand side of (\ref{q5}). By
construction, the operator $(\widehat{L}^*_\vartheta,
\mathcal{D}_\vartheta)$ is closed. Let $\mathcal{Q}_\vartheta$ be
the closure of $\mathcal{D}_\vartheta $ in $\mathcal{K}_\vartheta$.
Note that $\mathcal{Q}_\vartheta$ is a proper subset of
$\mathcal{K}_\vartheta$. By $L^\odot_\vartheta$ we define the part
of  $\widehat{L}^*_\vartheta$ in $\mathcal{Q}_\vartheta$. That is,
it is the restriction of $\widehat{L}^*_\vartheta$ to the set
\begin{equation}
  \label{q7}
  \mathcal{D}_\vartheta^\odot := \{ u \in \mathcal{D}_\vartheta: \widehat{L}^*_\vartheta u \in \mathcal{Q}_\vartheta\}.
\end{equation}
\begin{lemma}
  \label{1lm}
The operator $(L^\odot_\vartheta, \mathcal{D}_\vartheta^\odot)$ is the generator of a $C_0$-semigroup of bounded linear operators $S^\odot_\vartheta(t):
\mathcal{Q}_\vartheta \to \mathcal{Q}_\vartheta$, $t\geq 0$. Furthermore, for each $\vartheta' > \vartheta$, it follows that $\mathcal{K}_{\vartheta'} \subset \mathcal{D}_\vartheta^\odot$.
\end{lemma}
The proof of the lemma is given in the next section. Now we turn to
the problem in (\ref{12}). By the very definition of $L^\odot_\vartheta$, its action
 on the elements of $\mathcal{D}_\vartheta^\odot$
is given by the right-hand side of (\ref{q5}). Then the version of
(\ref{12}) in $\mathcal{Q}_\vartheta \subset \mathcal{G}_\vartheta$
is 
\begin{equation}
  \label{12a}
\frac{d}{dt} k_t = L^\odot_\vartheta k_t, \qquad k_t|_{t=0} = k_0
\in \mathcal{D}_\vartheta^\odot.
\end{equation}
\begin{definition}
  \label{2df}
By the classical global solution of the problem in (\ref{12a}) we
mean a function $[0,+\infty) \ni t \mapsto k_t \in
\mathcal{Q}_\vartheta$ which is continuously differentiable on
$[0,+\infty)$, lies in $\mathcal{D}_\vartheta^\odot$, and satisfies
(\ref{12a}).
\end{definition}
\begin{theorem}
  \label{1tm}
For each $\vartheta \in \mathbb{R}$  and $k_0 \in
\mathcal{D}_\vartheta^\odot$, the problem in (\ref{12}) has a unique global
classical solution $k_t\in \mathcal{Q}_\vartheta \subset
\mathcal{K}_\vartheta$ given by the formula
\begin{equation*}
k_t = S^{\odot}_\vartheta (t) k_0, \qquad t >0 ,
\end{equation*}
where $S^{\odot}_\vartheta$ is as in Lemma \ref{1lm}.  This, in
particular, holds if $k_0 \in \mathcal{K}_{\vartheta'}$ for some
$\vartheta' > \vartheta$. Furthermore, if $k_0$ is the correlation
function of some  $\mu_0\in \mathcal{P}(\Gamma)$, then, for each
$t>0$, there exists a unique $\mu_t\in \mathcal{P}(\Gamma)$ such
that $k_t$ is the the correlation function of this  $\mu_t$.
\end{theorem}
Note that in \cite{Berns} where the jumps were not free, the
evolution $\mathcal{K}_{\vartheta'} \ni k_0\mapsto k_t\in
\mathcal{K}_{\vartheta}$ was obtained on a bounded time interval
$[0, T)$ with $T$ dependent on the difference $\vartheta' -
\vartheta$.

\section{The Proofs}

\label{S3}

We first prove Lemma \ref{1lm} by means of a statement, which we
present here. Let $X$ be a Banach space with a cone of positive
elements, $X^+$, which is convex, generating ($X = X^+ - X^+$), and
proper ($X^+ \cap (-X^+)= \{0\}$). Let also the norm of $X$ be
additive on $X^+$, that is, $\|x+x'\|_X = \|x\|_X + \|x'\|_X$ for
$x, x' \in X^+$. Then there exists a positive (hence bounded) linear
functional $\varphi$ on $X$ such that $\varphi(x) =\|x\|_{X}$ for
each $x\in X^+$. Let now $X_1\subset X$ be a dense linear subset
equipped with its own norm in which it is also a Banach space, and
let the embedding $X_1 \hookrightarrow X$ be continuous. Assume also that
the norm of $X_1$ is additive on the cone $X_1^+:= X_1 \cap X^+$,
and $\varphi_1$ is the functional with property $\varphi_1(x)
=\|x\|_{X_1}$ for each $x\in X^+_1$. Let $S:= \{S(t)\}_{t\geq 0}$ be
a $C_0$-semigroup of bounded linear operators $S(t):X \to X$. It is
called sub-stochastic (resp. stochastic) if $S(t):X^+ \to X^+$ and
$\|S(t)\|_X \leq 1$ (resp. $\|S(t)\|_X = 1$) for all $t> 0$. Suppose
now that $(A_0, D(A_0))$ be the generator of a sub-stochastic
semigroup $S_0$ on $X$. Set $\check{S}_0 (t) = S_0(t)|_{X_1}$,
$t>0$, and assume that the following holds: \vskip.1cm
\begin{itemize}
  \item[(a)] for each $t>0$, $S_0(t) X_1 \to X_1$;
  \item[(b)] $\check{S}_0 := \{\check{S}_0 (t)\}_{t\geq 0}$ is a $C_0$-semigroup on $X_1$.
\end{itemize}
\vskip.1cm
Under these conditions, $\check{A}_0$ -- the generator of $\check{S}_0$, is the part of $A_0$ in $X_1$, see \cite[Proposition II.2.3]{EN}.
That is, $\check{A}_0$ is the restriction of $A_0$ to, cf. (\ref{q7}),
\begin{equation*}
D(\check{A}_0):= \{ x\in D(A_0)\cap X_1: A_0 x \in X_1\}.
\end{equation*}
The next statement is an adaptation of Proposition
2.6 and Theorem 2.7 of \cite{TV}.
\begin{proposition}
  \label{2pn}
Let conditions (a) and (b) given above hold, and $-A_0$ be a positive
linear operator in X. Let also $B$ be positive and such that its domain in $X$ contains $D(A_0)$ and
\begin{equation}
  \label{q9}
  \varphi((A_0 + B)x) = 0, \qquad x\in D(A_0)\cap X^+.
\end{equation}
Additionally, suppose that
$B: D(\check{A}_0) \to X_1$ and the following holds
\begin{equation}
  \label{q10}
 \varphi_1 ( (A_0 + B)x) \leq C \varphi_1 (x) - \varepsilon \|A_0 x\|_X , \qquad x \in D(\check{A}_0) \cap X^+,
\end{equation}
for some positive constants $C$ and $\varepsilon$. Then $(A,D(A))$
-- the closure of $(A_0 + B, D(A_0))$, is the generator of a
stochastic semigroup $S$ on $X$, which leaves $X_1$ invariant. That is, for each $t$, $S(t):X_1 \to X_1$.
\end{proposition}
\vskip.1cm \noindent
\subsection{Proof of Lemma \ref{1lm}}
\begin{proof}
 The space $\mathcal{G}_\vartheta$ possesses all the properties of the space $X$ assumed above. The corresponding functional is, cf. (\ref{q2a}),
\begin{equation}
  \label{q11}
  \varphi(v) = \int_{\Gamma_0} v (\eta) \exp(- \vartheta |\eta|) \lambda (d \eta).
\end{equation}
Then $(\widehat{A}_\vartheta, \widehat{\mathcal{G}}_\vartheta)$, see
(\ref{q3}) and (\ref{q3a}), generates  the sub-stochastic
$C_0$-semigroup $S_0$ of multiplication operators defined by the formula
\[
(S_0 (t) v)(\eta) = \exp( - \alpha |\eta|) v (\eta).
\]
Now let $\beta: \mathbb{N}_0 \to [0, +\infty)$ be such that $\beta (n) \to +\infty$ as $n\to +\infty$. Set
\begin{gather*}
\|v\|_{\vartheta,\beta } = \int_{\Gamma_0} |v(\vartheta)| \beta(|\eta|) \exp(- \vartheta |\eta|) \lambda (d \eta) \\[.2cm]
\varphi_\beta (v) = \int_{\Gamma_0}  v(\vartheta) \beta(|\eta|) \exp(- \vartheta |\eta|) \lambda (d \eta) , \nonumber \\[.2cm]
\mathcal{G}_{\vartheta,\beta}  = \{ v \in \mathcal{G}_{\vartheta,1}:
\|v\|_{\vartheta,\beta } < +\infty\} . \nonumber
\end{gather*}
Then clearly $S_0:\mathcal{G}_{\vartheta,\beta} \to
\mathcal{G}_{\vartheta,\beta}$, and $\| S_0 (t) v -
v\|_{\vartheta,\beta}  \to 0$ as $t\downarrow 0$ by the dominated
convergence theorem. Thus, both conditions (a) and (b) above are
satisfied. Next, let $\widehat{B}_\vartheta$ be as in (\ref{q3a}),
(\ref{l}). Then $(\widehat{L}_\vartheta,
\widehat{\mathcal{G}}_\vartheta)$, see (\ref{q5}), is closed.  As in
(\ref{q4}), by  (\ref{21}) and (\ref{q11}) we  obtain
\begin{eqnarray*}
\varphi(\widehat{L}_\vartheta v) & = & - \alpha \int_{\Gamma_0} |\eta | v(\eta)  \exp(- \vartheta |\eta|) \lambda (d \eta)\\[.2cm] & + & \int_{\Gamma_0} \int_{\mathbb{R}^d} \sum_{x\in \eta} a(x-y) v(\eta \setminus x \cup y)  \exp(-\vartheta |\eta|) \lambda ( d\eta)dy \nonumber \\[.2cm]
& =  & 0. \nonumber
\end{eqnarray*}
That is, the condition in (\ref{q9}) holds in our case. Moreover, again by (\ref{21})
we obtain that
\[
\varphi_\beta(\widehat{L}_\vartheta v) =0,
\]
hence, in our case (\ref{q10}) is satisfied if, for some $C>0$, the following holds
\begin{equation*}
 \beta (n) \geq C n, \qquad n \in \mathbb{N}.
\end{equation*}
Then we choose, e.g., $\beta (n) = n$ and obtain by Proposition \ref{2pn} that $(\widehat{L}_\vartheta, \widehat{\mathcal{G}}_\vartheta)$ generates a stochastic $C_0$-semigroup on
$\mathcal{G}_\vartheta$, which we denote by $\widehat{S}_\vartheta$. By the construction performed in (\ref{q6}) and (\ref{q7}) and by \cite[Theorem 10.4, page 39]{Pazy}, the semigroup in question is obtained as the restriction of the adjoint semigroup
$\widehat{S}^*_\vartheta$ to $\mathcal{Q}_\vartheta$.
\end{proof}

\subsection{Proof of Theorem \ref{1tm}}
\begin{proof}
The proof of the first part Theorem \ref{1tm} readily follows by
Lemma \ref{1lm} and \cite[Theorem 1.3, page 102]{Pazy}. So, it remains to prove that, for each $t>0$,
the solution $k_t$ is the correlation function for a unique $\mu_t\in \mathcal{P}(\Gamma)$. According to Proposition
\ref{1pn}, to this end we have to show that, for each $G\in B_{\rm bs}^+(\Gamma_0)$ and all $t>0$,
\begin{equation}
  \label{q16}
\langle \! \langle G, k_t \rangle \! \rangle = \langle \! \langle G,S^\odot_\vartheta (t) k_0 \rangle \! \rangle\geq 0,
\end{equation}
whenever this property holds for $t=0$. The proof of (\ref{q16}) follows along the following line of arguments.
For a $\Lambda \in \mathcal{B}_{\rm b}(\mathbb{R}^d)$, we let $$\Gamma_\Lambda = \{ \gamma \in \Gamma: \gamma \subset \Lambda\},$$ which is a Borel subset of $\Gamma$ such that $\Gamma_\Lambda \subset \Gamma_0$, see (\ref{1}) and (\ref{6}). Hence, $\Gamma_\Lambda \in \mathcal{B}(\Gamma_0)$. Next,  by $\mu_0^\Lambda$ we denote the projection of $\mu_0$ on $\Gamma_\Lambda$, i.e., $\mu^\Lambda_0   :=  \mu_0 \circ p_\Lambda^{-1}$, where $p_\Lambda (\gamma) = \gamma \cap \Lambda$ is the projection of $\Gamma$ onto $\Gamma_\Lambda$. It can be shown that $\mu^\Lambda_0$ is absolutely continuous with respect to the Lebesgue-Poisson measure $\lambda$. Let $R_0^\Lambda$ be its Radon-Nikodym derivative. For $N\in \mathbb{N}$, we also set
$R_0^{N,\Lambda} (\eta) = R_0^\Lambda (\eta) I_N (\eta)$, where $I_N (\eta) = 1$ if $|\eta|\leq N$ and $I_N (\eta) = 0$ otherwise. By this construction, $R_0^{N,\Lambda} \in \mathcal{G}_\theta$ for each $\theta \in \mathbb{R}$, see (\ref{q2b}). Thus, we can get
\begin{equation}
  \label{q17}
  R_t^{N,\Lambda} = \widehat{S}_\theta (t) R_0^{N,\Lambda}, \qquad t>0,
\end{equation}
where $\widehat{S}_\theta$ is the semigroup constructed in the proof of Lemma \ref{1lm}. Then $R_t^{N,\Lambda}(\eta) \geq 0$ and $\|R_t^{N,\Lambda}\|_{\theta, 1} \leq 1$.
Set
\begin{equation}
  \label{q18}
  q_t^{N,\Lambda} (\eta) = \int_{\Gamma_0} R_t^{N,\Lambda}(\eta \cup \xi) \lambda ( d\xi).
\end{equation}
By (\ref{bsa}), (\ref{bsb}),  and (\ref{21a}) for $G\in B^+_{\rm
bs}(\Gamma_0)$, we then get
\begin{equation}
  \label{q19}
 \langle \! \langle G , q_t^{N,\Lambda}  \rangle \! \rangle =  \langle \! \langle K G , R_t^{N,\Lambda} \rangle \! \rangle \geq 0.
\end{equation}
On  the other hand, by (\ref{q18}) for $t=0$ we have
\begin{equation*}
  0 \leq q_0^{N,\Lambda} (\eta) \leq \int_{\Gamma_0} R_t^{\Lambda}(\eta \cup \xi) \lambda ( d\xi) = k_{\mu_0}(\eta) \mathbb{I}_{\Gamma_\Lambda} ( \eta) \leq
   k_{\mu_0}(\eta),
\end{equation*}
where $\mathbb{I}_{\Gamma_\Lambda}$ is the corresponding indicator function.
Hence,  $q_0^{N,\Lambda} \in \mathcal{D}^\odot_\vartheta$ and we get
\begin{equation}
  \label{q21}
  k_t^{N,\Lambda} = S^\odot_\vartheta (t) q_0^{N,\Lambda}, \qquad t>0.
\end{equation}
As in \cite[Appendix]{Berns}, one can show that
\begin{equation}
  \label{q22}
\lim_{l\to +\infty} \lim_{n\to + \infty} \langle \! \langle G , k_t^{N_l,\Lambda_n}  \rangle \! \rangle =  \langle \! \langle G , k_t  \rangle \! \rangle ,
\end{equation}
for certain increasing sequences $\{N_l\}_{l\in \mathbb{N}}$ and $\{\Lambda_n\}_{n\in \mathbb{N}}$ such that $N_l\to +\infty$ and $\Lambda_n \to \mathbb{R}^d$. In (\ref{q22}), $k_t$ is the same as in (\ref{q16}).
Thus, by (\ref{q19}) and (\ref{q22}), we can obtain (\ref{q16}) by proving that
\begin{equation}
  \label{q23}
  \langle \! \langle G , q_t^{N,\Lambda}  \rangle \! \rangle =  \langle \! \langle G , k_t^{N,\Lambda}  \rangle \! \rangle , \qquad t>0.
\end{equation}
Set
\begin{equation*}
\phi(t) =    \langle \! \langle G , q_t^{N,\Lambda}  \rangle \! \rangle , \qquad \psi(t)=   \langle \! \langle G , k_t^{N,\Lambda}  \rangle \! \rangle. \end{equation*}
By (\ref{q21}) as well as by (\ref{12a}) and (\ref{q6}) we then get
\begin{equation}
  \label{q25}
\psi' (t) = \langle \! \langle \widehat{L}_\vartheta G , k_t^{N,\Lambda}  \rangle \! \rangle ,
\end{equation}
which makes sense since elements of $B_{\rm bs}(\Gamma_0)$ clearly belong to the domain
of $\widehat{L}_\vartheta$ for each $\vartheta \in \mathbb{R}$. Moreover, the operator $ \widehat{L}_\vartheta$, see (\ref{q5}), can be defined as a bounded operators acting from $\mathcal{G}_{\vartheta'}$ to $\mathcal{G}_{\vartheta}$, for each $\vartheta' < \vartheta$. Hence, one can define also its powers $ \widehat{L}^n_\vartheta: \mathcal{G}_{\vartheta'}\to \mathcal{G}_{\vartheta} $ such that, see \cite[eq. (3.21), page 1039]{Berns},
\begin{equation}
 \label{q26}
\|\widehat{L}^n_\vartheta\|_{\vartheta'\vartheta} \leq \left(\frac{2 \alpha n}{e (\vartheta - \vartheta')} \right)^n,
\end{equation}
 where $\| \cdot \|_{\vartheta'\vartheta}$ is the corresponding operator norm. Since
$B_{\rm bs}(\Gamma_0)\subset \mathcal{G}_{\vartheta'}$ for each $\vartheta'$, by the latter estimate we conclude that
$\psi$ can be continued to a function analytic in some neighborhood of the point $t=0$. Its derivatives are
\begin{equation}
  \label{q27}
  \psi^{(n)} (0) = \langle \! \langle \widehat{L}^n_\vartheta G , q_0^{N,\Lambda}  \rangle \! \rangle, \qquad n \in \mathbb{N}.
\end{equation}
On the other hand, by (\ref{q17}) and (\ref{q18}), as well as by the fact that the action of all operators like $\widehat{L}_\vartheta$ on their domains is the same as that in (\ref{l}), we obtain
\begin{equation*}
\phi' (t) = \langle \! \langle \widehat{L}_\vartheta G , q_t^{N,\Lambda}  \rangle \! \rangle ,
\end{equation*}
cf. (\ref{q25}) and (\ref{q26}), which yields that also $\phi$ can be continued to a function analytic in some neighborhood of the point $t=0$, with the derivatives, cf (\ref{q27}),
$ \phi^{(n)} (0) =  \psi^{(n)} (0)$. These facts readily yield (\ref{q23}), which completes the proof.

\end{proof}


\end{document}